\documentclass[11pt]{article}
\usepackage{amsmath}
\usepackage{amssymb}
\usepackage{amsthm}

\usepackage[normalem]{ulem}

\def\fpd#1#2{{\displaystyle\frac{\partial #1}{\partial #2}}}

\def\del{\nabla}
\def\cinfty#1{C^{\scriptscriptstyle\infty}(#1)}
\def\DV#1{{\rm D}^{\scriptscriptstyle V}_{#1}}
\def\DH#1{{\rm D}^{\scriptscriptstyle H}_{#1}}

\def\vectorfields#1{{\cal X}(#1)}

\def\tvectorfields{\vectorfields{\tau}}

\def\hook{\mathop{\hbox to 6pt{\hrulefill}
                      \hbox{\vrule\phantom{\vbox to 7pt{}}}}}

\def\DV#1{{\rm D}^{\scriptscriptstyle V}_{#1}}
\def\DH#1{{\rm D}^{\scriptscriptstyle H}_{#1}}
\def\CV#1{C^{\scriptscriptstyle V}_{#1}}
\def\G{\Gamma}
\def\tt{{\bf t}}
\def\H#1{{#1}^{\scriptscriptstyle H}}
\def\V#1{{#1}^{\scriptscriptstyle V}}

\def\R{\mathbb{R}}

\def\onehalf{{\textstyle\frac12}}

\def\onefourth{{\textstyle\frac14}}

\def\sode{{\sc sode}}

\usepackage{color}

\begin{document}

\title{Compatibility aspects of the method of phase synchronization for decoupling linear second-order differential equations}
\author{W.\ Sarlet$^\dagger$ and T.\ Mestdag$^{\dagger\,\ddagger}$\\[2mm]
{\small $^\dagger$ Department of Mathematics, Ghent University }\\
{\small Krijgslaan 281, 9000 Ghent, Belgium}\\[2mm]
{\small $^\ddagger$ Department of Mathematics,  University of Antwerp,}\\
{\small Middelheimlaan 1, 2020 Antwerpen, Belgium}
\\[2mm]
{\small Email: willy.sarlet@ugent.be, tom.mestdag@uantwerpen.be}
}

\date{\em Dedicated to Professor Tony Bloch on the occasion of his 65th birthday}

\maketitle

\begin{quote}
{\small {\bf Abstract.}
}
The so-called method of phase synchronization has been advocated in a number of papers as a way of decoupling a system of linear second-order differential equations by a linear transformation of coordinates and velocities. This is a rather unusual approach because velocity-dependent transformations in general do not preserve the second-order character of differential equations. Moreover, at least in the case of linear transformations, such a velocity-dependent one defines by itself a second-order system, which need not have anything to do, in principle, with the given system or its reformulation.
This aspect, and the related questions of compatibility it raises, seem to have been overlooked in the existing literature.
The purpose of this paper is to clarify this issue and to suggest topics for further research in conjunction with the general theory of decoupling in a differential geometric context.
\\[3mm]

\end{quote}

\section{Introduction}

In a series of papers, a number of authors have discussed techniques for decoupling the equations of motion of $n$-degree-of-freedom damped linear systems. We will take \cite{SKML18} as our starting reference, merely because it contains an accessible survey of the main formulas involved in the so-called method of phase synchronization, to which we can refer also for notations. But the original ideas which lead to this method go back to \cite{MIM09, MMI10, MM11, MMP11} and possibly more references therein. There is no doubt that the elegant explicit formulas concerning this phase synchronization have great merits for actually solving linear second-order equations, more particularly for developing efficient computer algorithms for doing so numerically. Yet, from an intrinsic mathematical point of view, there are two features about this literature which look puzzling. One is that there is a fairly complete looking literature about decoupling general (non-linear) second-order equations, which seems to have remained unnoticed by the authors in this field (see e.g.\ the review paper \cite{Sreview} for a survey of some  foundational references on this topic). Admittedly, this by itself is rather a side issue, because the method of phase synchronization for linear systems is a technique which would in a way be excluded from the start in this other part of the literature, where differential geometry is the proper environment for the study of second-order differential equations (\sode s). Briefly, a \sode\ lives on the tangent bundle of a differentiable manifold when it is autonomous, or on the first-jet bundle of a manifold fibred over $\R$ when explicit time-dependence is added to the picture. In both environments, it is natural to restrict coordinate transformations to be only of point-type, i.e.\ not depending on velocities, because these are the transformations which preserve the second-order character of \sode s. The linear transformations needed for phase synchronization, on the other hand, essentially depend on velocities. Yet, the authors in that field claim that the original linear \sode\ in such a way ``is transformed into a new one which is completely decoupled''. This brings us to the second puzzling feature referred to above: it cannot be quite the right terminology to say this! The point is this: there are in fact two \sode s which are associated to the transformation under consideration (under appropriate regularity assumptions), one is solely determined by the transformation in itself, the other one comes from applying the transformation to the given \sode\ indeed. For this reason, the whole procedure is only unambiguously defined when certain compatibility conditions are satisfied. That again is a feature which seems to have remained completely unnoticed in the literature.

We start in the next section with a brief summary of the formulas involved in the method of phase synchronization for solving linear \sode s. In Section~3, we analyze the compatibility problem and illustrate different aspects of it by a number of examples in Section~4.  In Section~5, we first present a brief overview of the differential geometric theory of separability of a \sode\ through an appropriate coordinate transformation. As indicated above, however, this cannot possibly cover what is happening in the method of phase synchronization. But we do succeed in presenting also a somewhat different differential geometric setting where such a mapping does make sense.

\section{Phase synchronization for linear systems}

To the best of our knowledge, the term {\em phase synchronization\/} was linked to a form of decoupling of linear \sode s for the first time in \cite{MIM09}. That may sound like a strange terminology for a method of decoupling, but its origin somehow simply comes from analyzing the solutions of the system under consideration. For the overlooked aspects of the method we wish to discuss in the next section, it is sufficient that we limit ourselves to the generic case of linear homogeneous \sode s with constant (real) coefficients. Also the case of so-called {\em defective systems\/} \cite{KMM11} is more like a side-issue here and will not be considered.

So, in the notations of \cite{SKML18}, consider an autonomous linear system
\begin{equation}
\mathbf{\ddot{q}} + C\,\mathbf{\dot{q}} + K\, \mathbf{q} = 0, \qquad \mathbf{q} \in \R^n. \label{1}
\end{equation}
Differential equations of this type model mechanical systems with both conservative ($-K\, \mathbf{q}$) and non-conservative forces ($- C\,\mathbf{\dot{q}}$). In the papers on the method of phase synchronization, it is often assumed that
 $C$  and $K$ are both symmetric; and the matrices are referred to as the damping and stiffness matrices. The classical terminology
is that we have dissipation of Rayleigh type (see e.g.\
\cite{Goldstein}). But, since systems (\ref{1}) with a gyroscopic force (where $C$ is skewsymmetric) are of equal importance in the literature (see e.g.\ \cite{Bloch,Rosenberg}), we will make  no assumption on $C$ (or $K$) here.

In the engineering literature, {\em modal analysis\/} is the term which is used when decoupling is possible by a point transformation. It requires that $C$ and $K$ commute plus some other conditions (cfr.\ our brief review in Section~5). Leaving that well-known situation aside, the basic idea of the phase synchronization approach is simply: start by looking for a particular solution of the form $\mathbf{q} = (\exp{\lambda t})\, \mathbf{v}$. This gives rise to the following quadratic eigenvalue problem:
\begin{equation}
(\lambda^2 \, I + \lambda\, C + K)\,\mathbf{v} = 0, \label{eigenv}
\end{equation}
and we assume for simplicity that all its solutions are distinct (real or complex). The claim in the technique of phase synchronization then is that there is a linear invertible transformation of the form
\begin{eqnarray}
\mathbf{q}\!\! &=&\!\! T_1\,\mathbf{p} + T_2\, \mathbf{\dot{p}}, \nonumber \\[-3mm]
     &   &  \label{tf} \\[-3mm]
\mathbf{\dot{q}} \!\!&=&\!\! T_3\,\mathbf{p} + T_4\, \mathbf{\dot{p}}, \nonumber
\end{eqnarray}
where the $T_i$ are real $n\times n$ matrices, which will transform (\ref{1}) into a decoupled system:
\begin{equation}
\mathbf{\ddot{p}} + D\,\mathbf{\dot{p}} + B\, \mathbf{p} = 0, \label{eqp}
\end{equation}
i.e.\ the matrices $D$ and $B$ will be real and diagonal. In fact, quite remarkably, an explicit solution is known for the matrices $T_i$ that will do the job. Various qualitative aspects of solving the quadratic eigenvalue problem (\ref{eigenv}) and coming to the ``state space formulation of phase synchronization'' (\ref{tf}) were discussed in \cite{MIM09, MMI10, MM11, MMP11}. The explicit results can be summarized as follows.

The quadratic eigenvalue problem (\ref{eigenv}) has $2n$ solutions. Complex ones appear in pairs (say $2c$) of complex conjugate numbers. To fix the idea, denote the ones with positive imaginary part, in de- or increasing order by $\lambda_1, \ldots, \lambda_c$, and their conjugates by $\lambda_{n+j} = \overline{\lambda}_j, \ j= 1, \ldots, c$. Split the remaining $2r$ ($r+c=n$) real eigenvalues also in pairs $(\lambda_j, \lambda_{n+j})$ for $j= c+1, \ldots, n$ (preferably in a certain order).
Denote corresponding eigenvectors by $\mathbf{v}_j, \mathbf{v}_{n+j}, \ j= 1, \ldots, n$. Now introduce the following matrices:
\begin{equation}
\Lambda_1 = \textrm{diag}\,[\lambda_1, \ldots, \lambda_n] \qquad\qquad \Lambda_2 = \textrm{diag}\,[\lambda_{n+1}, \ldots, \lambda_{2n}], \label{Lambda}
\end{equation}
and
\begin{equation}
V_1 = [\mathbf{v}_1, \ldots, \mathbf{v}_n] \qquad\qquad V_2 = [\mathbf{v}_{n+1}, \ldots, \mathbf{v}_{2n}]. \label{V}
\end{equation}
Then the claim is that the original equations (\ref{1}) are transformed into the decoupled ones (\ref{eqp}), with
\begin{equation}
D = - (\Lambda_1 + \Lambda_2), \qquad\qquad B = \Lambda_1\,\Lambda_2, \label{DB}
\end{equation}
and the transformation matrices $T_i$ which do the job are given by
\begin{eqnarray}
T_1 \!\!&=&\!\! (V_1\,\Lambda_2 - V_2\,\Lambda_1)(\Lambda_2 - \Lambda_1)^{-1} \nonumber \\
T_2 \!\!&=&\!\! (V_2 - V_1)(\Lambda_2 - \Lambda_1)^{-1} \nonumber \\[-3mm]
    & &  \label{T} \\[-3mm]
T_3 \!\!&=&\!\! (V_1 - V_2)(\Lambda_1\Lambda_2)(\Lambda_2 - \Lambda_1)^{-1} \nonumber \\
T_4 \!\!&=&\!\! (V_2\,\Lambda_2 - V_1\,\Lambda_1)(\Lambda_2 - \Lambda_1)^{-1}. \nonumber
\end{eqnarray}
Explicit expressions for the matrices of the inverse transformation of (\ref{tf}) can be found e.g.\ in \cite{SKML18} as well.

In retrospect, it should perhaps not come as a surprise that every linear system like (\ref{1}) can be replaced by a set of decoupled equations, once you have managed to solve the original equations explicitly. But what is extra appealing here is the availability of an explicit so-called state space transformation (\ref{tf},\! \ref{T}), associated to this solution. That is why we think it is relevant to understand and be aware of some features of this transformation which have been hidden so far.

\section{Compatibility aspects}

Let's have a closer look at transformations of the form (\ref{tf}) with constant coefficient matrices in all generality. Taking a derivative of the first relation $\mathbf{q} = T_1\,\mathbf{p} + T_2\, \mathbf{\dot{p}}$, we get
\begin{equation}
\mathbf{\dot{q}} = T_1\,\mathbf{\dot{p}} + T_2\, \mathbf{\ddot{p}}. \label{dotq}
\end{equation}
Identification with the second relation $\mathbf{\dot{q}} = T_3\,\mathbf{p} + T_4\, \mathbf{\dot{p}}$ then leads to
\begin{equation}
T_2\,\mathbf{\ddot{p}} + (T_1 - T_4)\,\mathbf{\dot{p}} - T_3\, \mathbf{p} = 0. \label{sode1}
\end{equation}
Hence, provided that $T_2$ is non-singular, the transformation under consideration defines a regular \sode\ by itself. This is of course a specific feature of the non-point character of the transformation, since (\ref{sode1}) will not occur when $T_2=T_3=0$ and $T_4=T_1$.

On the other hand, taking a derivative of the second relation, we get
\begin{equation}
\mathbf{\ddot{q}} = T_3\,\mathbf{\dot{p}} + T_4\, \mathbf{\ddot{p}}. \label{ddotq}
\end{equation}
Inserting (\ref{tf}) and (\ref{ddotq}) into the original equation (\ref{1}), it transforms into
\begin{equation}
T_4\,\mathbf{\ddot{p}} + (T_3 + C\,T_4 + K\,T_2)\,\mathbf{\dot{p}} + (C\,T_3 + K\,T_1)\, \mathbf{p} = 0, \label{sode2}
\end{equation}
which gives another regular \sode\ provided that $T_4$ is non-singular, and could be termed the transformed equation.

For this to make any sense, however, it is obvious that (\ref{sode1}) and (\ref{sode2}) must be compatible, in one of the two possible interpretations which will now be explored. In the generic case that both $T_2$ and $T_4$ are non-singular, the compatibility conditions read:
\begin{eqnarray}
{T_2}^{-1} (T_1 - T_4) \!\!&=&\!\! {T_4}^{-1} (T_3 + C\,T_4 + K\,T_2)  \nonumber \\[-3mm]
   & &  \label{cc1} \\[-3mm]
- {T_2}^{-1} T_3 \!\!&=&\!\! {T_4}^{-1} (C\,T_3 + K\,T_1). \notag
\end{eqnarray}
So in general, there seems to be some freedom in this process since we have two conditions for four unknowns. It might be tempting, for example, to say: $T_2$ and $T_4$ can be chosen arbitrarily (non-singular) and then (\ref{cc1}) serves to find meaningful options for $T_1$ and $T_3$. Of course, we are not necessarily talking then about the specific purpose of arriving at decoupled equations.

Going back to the strategy leading to the specific transformation (\ref{T}), it can be verified, always under the same regularity assumptions, that these explicit expressions for the $T_i$ effectively satisfy the conditions (\ref{cc1}). However, the following immediate remarks naturally pop up then:

\begin{itemize}\vspace{-2mm}
\item[-] The freedom of choice of a factor in selecting eigenvectors always makes it possible to make either $T_2$ or $T_4$ singular.
\item[-] It is not difficult to construct an example where both $T_2$ and $T_4$ are singular!
\end{itemize}

These are elements which make that it is not quite right to say that a `state space transformation' like (\ref{tf}) {\em transforms\/} the original \sode\ (\ref{1}) into decoupled equations (\ref{eqp}). Let's think for a moment about the case $n=2$ for simplicity. When $T_4$ is singular, for example, the `transformed equation' (\ref{sode2}) contains only one expression and needs the assistance of (\ref{sode1}) to arrive at a 2-dimensional \sode\ again. And when both $T_2$ and $T_4$ are singular, both (\ref{sode1}) and (\ref{sode2}) have only one independent expression, each possibly involving $p_1$ and $p_2$. It is then only after joining these two that we see a system emerge that can be rewritten into a decoupled form (see the examples in the next section).

The above observations prompt us for a second possible way of looking at the compatibility conditions. In fact, the claim is essentially that there exist two diagonal matrices $D$ and $B$, say of the form (\ref{DB}) which define a decoupled system (\ref{eqp}) associated to the original \sode\ (\ref{1}). Substitution of (\ref{eqp}) into (\ref{sode1}) and (\ref{sode2}) leads to the requirements
\begin{eqnarray}
T_2\,\mathbf{\ddot{p}} + (T_1 - T_4)\,\mathbf{\dot{p}} - T_3\, \mathbf{p} \!\!&\equiv &\!\!  T_2 \, (\mathbf{\ddot{p}} + D\,\mathbf{\dot{p}} + B\, \mathbf{p}), \notag\\[-3mm]
   & &  \label{cc2} \\[-3mm]
T_4\,\mathbf{\ddot{p}} + (T_3 + C\,T_4 + K\,T_2)\,\mathbf{\dot{p}} + (C\,T_3 + K\,T_1)\, \mathbf{p} \!\!&\equiv &\!\!  T_4 \,(\mathbf{\ddot{p}} + D\,\mathbf{\dot{p}} + B\, \mathbf{p}). \notag
\end{eqnarray}
Identification of coefficients then requires
\begin{equation}
T_3 = - T_2 B \qquad \mbox{and} \qquad T_1 = T_4 + T_2 D \label{id1}
\end{equation}
from the first, and likewise from the second:
\begin{equation}
T_4 D = T_3 + C T_4 + K T_2 \quad \mbox{and} \quad
T_4 B = C T_3 + K T_1. \label{id2}
\end{equation}
From this perspective, the compatibility requirements (\ref{cc1}) are replaced, for example, by conditions which follow from substituting $T_3$ and $T_1$ from (\ref{id1}) into (\ref{id2}), which results in
\begin{eqnarray}
T_4 D + T_2 B \!\!&=&\!\! C T_4 + K T_2, \notag\\[-3mm]
   & &  \label{cc3} \\[-3mm]
T_4 B \!\!&=&\!\! - C T_2 B + K (T_4 + T_2 D). \notag
\end{eqnarray}
There is then no problem to deal with situations where $T_2$ and/or $T_4$ are singular, but the main point we are making here is that there still is an interesting compatibility issue. In all generality, it could be the main tool for finding suitable transformations, whatever the specific objective is for the matrices $D$ and $B$ in the new equation. But let us stick to the purpose of decoupling here, meaning that $D$ and $B$ must be diagonal and thus can just as well be represented as sum and product of as yet undetermined diagonal matrices $\Lambda_i$ as in (\ref{DB}). We shall illustrate that one can then obtain the quadratic eigenvalue problem (\ref{eigenv}) purely from manipulations of the conditions (\ref{cc3}).

Substituting (\ref{DB}) into (\ref{cc3}), we are looking at linear homogeneous equations for $T_2$ and $T_4$, which will have solutions provided some $2n\times 2n$ determinant vanishes. But look what is happening if we introduce `new variables', $W_1$ and $W_2$ say, by the following relations:
\begin{eqnarray}
T_2 \!&=&\! W_2 - W_1, \notag \\[-3mm]
    & & \label{W1W2} \\[-3mm]
T_4 \!&=&\! W_2 \Lambda_2 - W_1 \Lambda_1. \notag
\end{eqnarray}
Needless to say, the transition $(T_2,T_4) \leftrightarrow (W_1,W_2)$ is inspired by the results (\ref{T}) and we observe that it is invertible, provided that $\Lambda_2 - \Lambda_1$ is non-singular. The substitution at first gives rise to the equations
\begin{eqnarray}
W_1 {\Lambda_1}^2 - W_2 {\Lambda_2}^2 \!&=&\! C(W_2\Lambda_2 - W_1\Lambda_1) + K(W_2 - W_1), \notag \\[-3mm]
    & & \label{cc4} \\[-3mm]
W_2\Lambda_1{\Lambda_2}^2 - W_1{\Lambda_1}^2\Lambda_2 \!&=&\! C(W_1-W_2)\Lambda_1\Lambda_2 + K(W_1\Lambda_2 - W_2\Lambda_1). \notag
\end{eqnarray}
However, multiplying the first of these on the right by $\Lambda_1$ (respectively $\Lambda_2$) and adding the second, we obtain an equivalent system which reveals a rather unexpected effect! Indeed, both equations now look the same. That is to say, one is obtained from the other by simply interchanging the indices 1 and 2. Explicitly, after right multiplication by $(\Lambda_2 - \Lambda_1)^{-1}$, we in fact need two different solutions $(W_1,\Lambda_1)$ and $(W_2,\Lambda_2)$ of the single equation
\begin{equation}
W\Lambda^2 + CW\Lambda + KW = 0. \label{W}
\end{equation}
This is a drastic simplification, because the restrictions on suitable diagonal matrices $\Lambda$ now come from a vanishing $n\times n$ determinant. And if we look at what (\ref{W}) tells us for each column of $W$, we see that we are in fact facing the quadratic eigenvalue problem (\ref{eigenv}).

As is already clear from the preceding considerations (and of course also for example in \cite{MMI10}), for a given \sode\ (\ref{1}) and its decoupled form (\ref{eqp}), a transformation (\ref{tf}) which does the job is by far not unique. What's more, due to the freedom of choice we already have when selecting the matrices of eigenvalues $\Lambda_1$ and $\Lambda_2$, one and the same \sode\ (\ref{1}) can even be transformed into different decoupled equations of the form (\ref{eqp}).
In the next section we present a number of examples which illustrate most of the points we have raised so far.

\section{Illustrative examples}

{\bf Example 1.}\ Let's start by picking an example which has been discussed in \cite{MMI10}. Suppose we have a 2-dimensional system (\ref{1}) with
\begin{equation}
C = \begin{pmatrix} 4 & -1 \\ -1 & 8 \end{pmatrix} \qquad K = \begin{pmatrix} 1 & 0 \\ 0 & 4 \end{pmatrix}. \label{ex1}
\end{equation}
The solutions of the quadratic eigenvalue problem (\ref{eigenv}) are real and can be used to define matrices $\Lambda_1$ and $\Lambda_2$ for a decoupled system, for example as follows:
\begin{equation}
 \Lambda_1= \begin{pmatrix} -2+\sqrt{2} & 0 \\ 0 & -4+\sqrt{14} \end{pmatrix} \qquad \Lambda_2 = \begin{pmatrix} -2-\sqrt{2} & 0 \\ 0 & -4-\sqrt{14} \end{pmatrix}. \label{ex1-1}
\end{equation}
Following the pattern of the relations (\ref{V}\,-\,\ref{T}), we select corresponding eigenvectors and put
\begin{equation}
 V_1= \begin{pmatrix} -2+\sqrt{2} & -4+\sqrt{14} \\ -1 & 15-4\sqrt{14} \end{pmatrix} \qquad V_2 = \begin{pmatrix} -2-\sqrt{2} & -4-\sqrt{14} \\ -1 & 15+4\sqrt{14} \end{pmatrix}. \label{Vi1-1}
\end{equation}
This results in transformation matrices $T_i$, given by
\begin{equation}
 T_1= \begin{pmatrix} 0 & 0 \\ -1 & -1 \end{pmatrix} \quad T_2 = \begin{pmatrix} 1 & 1 \\ 0 & -4 \end{pmatrix} \quad T_3 = \begin{pmatrix} -2 & -2 \\ 0 & 8 \end{pmatrix} \quad T_4 = \begin{pmatrix} -4 & -8 \\ -1 & 31 \end{pmatrix}. \label{Ti1-1}
\end{equation}
The main point to observe now is that both $T_2$ and $T_4$ are non-singular here. Hence, the transformation by itself defines a regular \sode\ (\ref{sode1}) and there is another one (\ref{sode2}) which comes from the given equations. We leave it to the reader to verify that the two compatibility conditions (\ref{cc1}) are effectively satisfied.

The eigenvectors in the columns of the $V_i$ in (\ref{Vi1-1}) were chosen by sort of focussing on the first line in the vanishing determinant for the eigenvalues. If we rather look at the second line, another natural choice for the $V_i$ reads:
\begin{equation}
 V_1= \begin{pmatrix} -6+4\sqrt{2} & 2 \\ -2+\sqrt{2} & -4+\sqrt{14} \end{pmatrix} \qquad V_2 = \begin{pmatrix} -6-4\sqrt{2} & 2 \\ -2-\sqrt{2} & -4-\sqrt{14} \end{pmatrix}. \label{Vi1-2}
\end{equation}
Corresponding $T_i$ which equally satisfy the compatibility conditions (\ref{cc1}) then become:
\begin{equation}
 T_1= \begin{pmatrix} 2 & 2 \\ 0 & 0 \end{pmatrix} \quad T_2 = \begin{pmatrix} 4 & 0 \\ 1 & 1 \end{pmatrix} \quad T_3 = \begin{pmatrix} -8 & 0 \\ -2 & -2 \end{pmatrix} \quad T_4 = \begin{pmatrix} -14 & 2 \\ -4 & -8 \end{pmatrix}. \label{Ti1-2}
\end{equation}
Choosing for a mixture of the two ideas has the special effect here that $T_2$ becomes the unit matrix. We then have
\begin{equation}
 V_1= \begin{pmatrix} -2+\sqrt{2} & 2 \\ -1 & -4+\sqrt{14} \end{pmatrix} \qquad V_2 = \begin{pmatrix} -2-\sqrt{2} & 2 \\ -1 & -4-\sqrt{14} \end{pmatrix}, \label{Vi1-3}
\end{equation}
and correspondingly
\begin{equation}
 T_1= \begin{pmatrix} 0 & 2 \\ -1 & 0 \end{pmatrix} \quad T_2 = \begin{pmatrix} 1 & 0 \\ 0 & 1 \end{pmatrix} \quad T_3 = \begin{pmatrix} -2 & 0 \\ 0 & -2 \end{pmatrix} \quad T_4 = \begin{pmatrix} -4 & 2 \\ -1 & -8 \end{pmatrix}. \label{Ti1-3}
\end{equation}
All three selections produce the same decoupled system (\ref{eqp}) with
\begin{equation}
 D = \begin{pmatrix} 4 & 0 \\ 0 & 8 \end{pmatrix} \qquad B = \begin{pmatrix} 2 & 0 \\ 0 & 2 \end{pmatrix}. \label{ex1-DB1}
\end{equation}

But let us also illustrate here that even the form of the decoupled equations depends on the freedom we have in the process. In this example, we might just as well combine the eigenvalues to define the $\Lambda_i$ as
\begin{equation}
 \Lambda_1= \begin{pmatrix} -2+\sqrt{2} & 0 \\ 0 & -2-\sqrt{2} \end{pmatrix} \qquad \Lambda_2 = \begin{pmatrix} -4+\sqrt{14} & 0 \\ 0 & -4-\sqrt{14} \end{pmatrix}. \label{ex1-2}
\end{equation}
The decoupled $\bf{p}$-equations then are determined by
\begin{eqnarray}
 D &=& \begin{pmatrix} 6 -\sqrt{2}-\sqrt{14} & 0 \\ 0 &  6 +\sqrt{2}+\sqrt{14}\end{pmatrix} \nonumber \\[-2mm]
   & &  \label{ex1-DB2} \\[-2mm] B &=& \begin{pmatrix} 8 - 2\sqrt{14} - 4\sqrt{2} + 2\sqrt{7} & 0 \\ 0 & 8 +2\sqrt{14} +4\sqrt{2} + 2\sqrt{7} \end{pmatrix}. \nonumber
\end{eqnarray}
A possible selection of associated eigenvectors leads to
\begin{equation}
 V_1= \begin{pmatrix} -2+\sqrt{2} & -2-\sqrt{2} \\ -1 & -1 \end{pmatrix} \qquad V_2 = \begin{pmatrix} -4+\sqrt{14} & -4-\sqrt{14} \\ 15-4\sqrt{14} & 15 +4\sqrt{14} \end{pmatrix}, \label{Vi1-4}
\end{equation}
and we spare the reader the exotic expressions of the resulting matrices $T_i$.

{\bf Example 2.}\ Consider the following \sode:
\begin{eqnarray}
\ddot{q}_1 \!\!&=& \!\!\dot{q}_1 + \dot{q}_2 - q_1 - \onehalf q_2 \nonumber \\[-3mm]
   & &  \label{ex2} \\[-3mm]
\ddot{q}_2 \!\!&=&\!\! \dot{q}_2 - 4\dot{q}_1 - q_2 + 2q_1. \nonumber
\end{eqnarray}
This is the linear subsystem of an example mentioned in \cite{ST, S} for reasons which will be briefly touched upon in the next section. We have
\begin{equation}
C = \begin{pmatrix} -1 & -1 \\ 4 & -1 \end{pmatrix} \qquad K = \begin{pmatrix} 1 & 1/2 \\ -2 & 1 \end{pmatrix}. \label{ex2CK}
\end{equation}
An interesting feature here is that $C$ and $K$ actually commute. Yet, `modal analysis' does not apply because $C$ is not diagonalizable over the reals (cfr. next section).

The solutions of the quadratic eigenvalue problem (\ref{eigenv}) are complex here. We need to make sure that the matrices $\Lambda_1$ and $\Lambda_2$ are complex conjugate, so that their sum and product in the resulting ${\bf p}$-equations (\ref{eqp}) become real. So, let's take
\begin{equation}
\Lambda_1 = \begin{pmatrix} \frac{1}{2} + \frac{1}{2}\, i(2 + \sqrt{7})  & 0 \\ 0 & \frac{1}{2} + \frac{1}{2}\, i(\sqrt{7}-2) \end{pmatrix} \qquad
\Lambda_2 = \overline{\Lambda_1}. \label{ex2-1}
\end{equation}
In the same way, if we arrange corresponding eigenvectors in such a way that the matrices $V_1$ and $V_2$ are complex conjugate, then the structure of the formulas (\ref{T}) will guarantee that also all $T_i$ become real. Let's choose
\begin{equation}
V_1= \begin{pmatrix} -\frac{1}{2} + \frac{1}{2}\,i & -\frac{1}{2} + \frac{1}{2}\,i \\ -1 - i & 1 + i \end{pmatrix} \qquad V_2 = \overline{V_1}. \label{Vi2-1}
\end{equation}
Then, the key matrices $T_2$ and $T_4$ for the \sode s (\ref{sode1}) and (\ref{sode2}) are non-singular and read:
\begin{equation}
T_2 = \begin{pmatrix} \frac{1}{3}(\sqrt{7} -2) & \frac{1}{3}(\sqrt{7} +2) \\[1mm] - \frac{2}{3}(\sqrt{7} -2) & \frac{2}{3}(\sqrt{7} +2) \end{pmatrix} \quad  T_4 = \begin{pmatrix} - \frac{1}{6}(5 -\sqrt{7}) & - \frac{1}{6}(1 -\sqrt{7})\\[1mm] - \frac{1}{3}(1 +\sqrt{7}) & \frac{1}{3}(5 +\sqrt{7}) \end{pmatrix}. \label{Ti2-1}
\end{equation}
For completeness, we further have
\begin{equation}
T_1 = \begin{pmatrix} -\frac{1}{6}(1+\sqrt{7}) & -\frac{1}{6}(5+\sqrt{7}) \\[1mm] - \frac{1}{3}(5-\sqrt{7}) & \frac{1}{3}(1-\sqrt{7}) \end{pmatrix} \quad  T_3 = \begin{pmatrix} - \frac{1}{3}(1+\sqrt{7}) &  \frac{1}{3}(1 -\sqrt{7})\\[1mm]  \frac{2}{3}(1 +\sqrt{7}) & \frac{2}{3}(1 -\sqrt{7}) \end{pmatrix}. \label{Ti2-2}
\end{equation}
It is tedious but straightforward to verify the compatibility conditions, and the decoupled ${\bf p}$-equations (\ref{eqp}) become
\begin{eqnarray}
\ddot{p}_1 \!\!&=&\!\! \dot{p}_1 - (3+\sqrt{7})\,p_1 \nonumber \\[-3mm]
   & &  \label{ex2-p} \\[-3mm]
\ddot{p}_2 \!\!&=&\!\! \dot{p}_2 - (3-\sqrt{7})\,p_2. \nonumber
\end{eqnarray}

Remark: It is not difficult here to choose eigenvectors in such a way that either $T_2$ or $T_4$ become singular, so that compatibility can no longer be checked via (\ref{cc1}). But we will give a more systematic discussion of such issues in the next example.

{\bf Example 3.}\ Consider a system (\ref{1}) with
\begin{equation}
C = \begin{pmatrix} 2 & 0 \\ -2 & 2 \end{pmatrix} \qquad K = \begin{pmatrix} 4 & 0 \\ -2 & 2 \end{pmatrix}. \label{ex3}
\end{equation}
All solutions of the quadratic eigenvalue problem are complex now. In such a case, the freedom of using them in choosing $\Lambda_i$ is of course not unlimited, as we wish sum and product of the $\Lambda_i$ to be real in accordance with (\ref{DB}). So let's take
\begin{equation}
 \Lambda_1= \begin{pmatrix} -1+i & 0 \\ 0 & -1+i\sqrt{3} \end{pmatrix} \qquad \Lambda_2 = \begin{pmatrix} -1-i & 0 \\ 0 & -1-i\sqrt{3} \end{pmatrix}. \label{ex3-1}
\end{equation}
But there remains a lot of freedom in selecting appropriate matrices $V_i$ in view of the free factor for each eigenvector. Let's put
\begin{equation}
 V_1= \begin{pmatrix} 0 & \mu_1 \\ \rho_1 & -i\sqrt{3}\,\mu_1 \end{pmatrix} \qquad V_2 = \begin{pmatrix} 0 & \mu_2 \\ \rho_2 & i\sqrt{3}\,\mu_2 \end{pmatrix}, \label{Vi3-1}
\end{equation}
with the $\rho_i, \mu_i$ as yet undetermined. The main issue next is to look at the structure of $T_2$ and $T_4$ in (\ref{T}). We want the transformation (\ref{tf}) to be real as well and if we make sure that this is the case for $T_2$ and $T_4$, it follows from (\ref{id1}) that $T_1$ and $T_3$ will be real also. It is fairly easy to see that for this purpose, it suffices to require that both the couple $(\rho_1,\rho_2)$ and the couple $(\mu_1,\mu_2)$ are equal if real, and complex conjugate when complex. We shall use this example to illustrate in greater detail the effects of having singularity in $T_2$ and/or $T_4$.

First we choose the free factors to be all real, for example simply $\rho_1=\rho_2=\mu_1=\mu_2=1$. Then it is easy to verify that the resulting $T_i$ matrices become:
\begin{equation}
 T_1= \begin{pmatrix} 0 & 1 \\ 1 & -1 \end{pmatrix} \quad T_2 = \begin{pmatrix} 0 & 0 \\ 0 & -1 \end{pmatrix} \quad T_3 = \begin{pmatrix} 0 & 0 \\ 0 & 4 \end{pmatrix} \quad T_4 = \begin{pmatrix} 0 & 1\\ 1 & 1 \end{pmatrix}. \label{Ti3-1}
\end{equation}
$T_2$ is singular while $T_4$ is not. Therefore, the equations (\ref{sode2}) define a regular \sode\ which should decouple when written in normal form; (\ref{sode1}) on the other hand will produce only one equation which is hopefully consistent with the other two. Explicitly, when we write the $\bf{p}$-equations exactly in the formats (\ref{sode1}), (\ref{sode2}), we get the single equation
\begin{equation}
- \ddot{p}_2 - 2\dot{p}_2 - 4 p_2 =0  \label{s1-31}
\end{equation}
for the first, and the system
\begin{eqnarray}
\ddot{p}_2 + 2\dot{p}_2 + 4 p_2 &=& 0   \nonumber \\[-3mm]
   & &  \label{s2-31} \\[-3mm]
\ddot{p}_1 + \ddot{p}_2 + 2(\dot{p}_1 + \dot{p}_2) + 2p_1 + 4p_2 &=& 0 \nonumber
\end{eqnarray}
for the second. The latter indeed produces the decoupled \sode\ (\ref{eqp}) and (\ref{s1-31}) is compatible as expected.

Secondly, let's choose factors $\rho_i$, $\mu_i$ which are all complex, say
\begin{equation}
\rho_1 = 1+i, \quad \rho_2=1-i, \quad \mu_1= 1+\sqrt{3}, \quad \mu_2= 1 - i\sqrt{3}. \label{factors2}
\end{equation}
Then, the resulting $V_i$ from (\ref{Vi3-1}) give rise to the following $T_i$:
\begin{equation}
 T_1= \begin{pmatrix} 0 & 2 \\ 2 & 2 \end{pmatrix} \quad T_2 = \begin{pmatrix} 0 & 1 \\ 1 & -1 \end{pmatrix} \quad T_3 = \begin{pmatrix} 0 & -4 \\ -2 & 4 \end{pmatrix} \quad T_4 = \begin{pmatrix} 0 & 0\\ 0 & 4 \end{pmatrix}. \label{Ti3-2}
\end{equation}
Clearly now, $T_4$ is singular while $T_2$ is not. This time, it is (\ref{sode1})
\begin{eqnarray}
\ddot{p}_2 + 2\dot{p}_2 + 4 p_2 &=& 0   \nonumber \\[-3mm]
   & &  \label{s1-32} \\[-3mm]
\ddot{p}_1 - \ddot{p}_2 + 2(\dot{p}_1 - \dot{p}_2) + 2p_1 - 4p_2 &=& 0 \nonumber
\end{eqnarray}
which generates a decoupled \sode, and (\ref{sode2}) is compatible again.

Finally, we try a mixture of complex and real proportionality factors, say $\rho_i = 1\pm i, \mu_i=1$. This time, we get
\begin{equation}
 T_1= \begin{pmatrix} 0 & 1 \\ 2 & -1 \end{pmatrix} \quad T_2 = \begin{pmatrix} 0 & 0 \\ 1 & -1 \end{pmatrix} \quad T_3 = \begin{pmatrix} 0 & 0 \\ -2 & 4 \end{pmatrix} \quad T_4 = \begin{pmatrix} 0 & 1\\ 0 & 1 \end{pmatrix}. \label{Ti3-3}
\end{equation}
So both $T_2$ and $T_4$ are singular now, meaning that neither (\ref{sode1}) nor (\ref{sode2}) are proper \sode s: (\ref{sode1}) gives us the second of (\ref{s1-32}), and (\ref{sode2}) the first and joining forces they produce the decoupled system.

\section{The differential geometric results about separability of \sode s}

The study of separability of \sode s in a differential geometric context is an interesting application of the theory of derivations of vector-valued forms along the tangent bundle projection, which was launched in the PhD work of Eduardo Mart\'{\i}nez. The interested reader might have a look at \cite{MCS1} and \cite{MCS2} for the basics of this theory. Necessary and sufficient conditions for separability of a \sode\ were first obtained in a strictly autonomous framework in \cite{MCS3}, and later extended to time-dependent situations in \cite{CSVM}. A key issue in this theory is the diagonalizability of a type (1,1) tensor field $\Phi$ and only real eigenfunctions were taken into consideration in the preceding references. A further extension of the theory for the case of complex eigenfunctions was developed in \cite{ST}.

As said in the introduction, there is no direct overlap in those theories with the content of the preceding sections. We will therefore limit ourselves here, for comparison, to a brief sketch of the theory in the autonomous set-up and see what it tells us about the special case of linear systems. But we do need to introduce some basic intrinsic operations in order to be able to formulate the relevant theorems.

An autonomous \sode\
\begin{equation}
\ddot{q}^i = f^i(q,\dot{q}) \qquad i=1,\ldots,n \label{nlsode}
\end{equation}
is governed by a vector field $\G$ on the tangent bundle $\tau:TM\rightarrow M$ (coordinates $(q,v)$) of a differentiable manifold $M$, which has the form (with summation convention)
\begin{equation}
\G = v^i \fpd{}{q^i} + f^i(q,v) \fpd{}{v^i}. \label{Gamma}
\end{equation}
$\G$ intrinsically defines a {\sl horizontal distribution\/}, spanned by the vector fields
\begin{equation}
H_i = \fpd{}{q^i} - \G^j_i\fpd{}{v^j}, \qquad \mbox{where} \quad \G^j_i = -\frac{1}{2}\fpd{f^j}{v^i}. \label{Hi}
\end{equation}
That horizontal distribution defines a so-called {\sl non-linear connection\/} on the tangent bundle and the $\G^j_i$ are its connection coefficients. In turn, there is an associated {\sl linear connection\/}
on what is called the pull-back bundle $\tau^*TM\rightarrow TM$, which basically has three constituents: a {\sl horizontal\/} and {\sl vertical covariant derivative\/}, denoted by $\DH{X}\ $ and $\DV{X}$ with $X \in\tvectorfields$, the set of vector fields along $\tau$, plus what is called the {\sl dynamical covariant derivative\/} $\del$. These act as degree zero derivations on arbitrary vector-valued forms along the tangent bundle projection $\tau$. Tensor fields along $\tau$ look like tensors on the base manifold $M$, but with coefficients which depend on the coordinates $(q,v)$ of the full space $TM$. $X \in\tvectorfields$, for example, has the coordinate representation
$X= X^i(q,v)\partial/\partial q^i$. The last thing we need for introducing all relevant objects is the concept of horizontal and vertical lifts from $\tvectorfields$ to $\vectorfields{TM}$. With $V_i$ as shorthand for the coordinate vectorfields $\partial/\partial v^i$, they are defined by
\begin{equation}
\H{X}=X^i\,H_i  \qquad \V{X}=X^i\,V_i=X^i \qquad \forall X \in \tvectorfields. \label{HV}
\end{equation}
$\DV{X}$ and $\DH{X}$  can now be defined by decomposition of the Lie-bracket of a horizontal and vertical lift into its horizontal and vertical part:
\begin{equation}
{}[\H{X},\V{Y}] = \V{(\DH{X}Y)} - \H{(\DV{Y}X)}. \label{DHDV}
\end{equation}
Likewise, $\del$ and $\Phi$ referred to above are fully determined by the bracket decomposition
\begin{equation}
{}[\G,\H{X}] = \H{(\del X)} + \V{\Phi(X)}. \label{delphi}
\end{equation}
Coordinate expressions will now further clarify these concepts. For the derivations, it is sufficient to know that
\begin{equation}
\DV{X}(F)=\V{X}(F), \quad \DH{X}=\H{X}(F), \quad \del(F)=\G(F), \qquad \mbox{with}\ F\in\cinfty{TM} \label{onF}
\end{equation}
and
\begin{equation}
\DV{X}\,\fpd{}{q^i}=0, \qquad \DH{X}\,\fpd{}{q^i}=
\V{X}(\G^j_i)\fpd{}{q^j}, \qquad \del\,\fpd{}{q^i}=\G^j_i\fpd{}{q^j}, \label{onvectors}
\end{equation}
while the action on 1-forms along $\tau$ subsequently follows by duality. As said before, an important element for our statements about the \sode\ $\G$ is the associated type (1,1) tensor field $\Phi$ along $\tau$ which pops up in (\ref{delphi}). It is called the {\sl Jacobi endomorphism\/} and has components
\begin{equation}
\Phi^i_j = -\fpd{f^i}{q^j} - \G^i_k\G^k_j - \G(\G^i_j). \label{Phi}
\end{equation}
$\Phi$ determines the curvature $R$ of the non-linear \sode-connection:
\begin{equation}
3\,R(X,Y) =  \DV{X}\Phi(Y) - \DV{Y}\Phi(X). \label{R}
\end{equation}
Another derived type (1,2) tensor of interest is
\begin{equation}
\CV{\Phi}(X,Y) = [\DV{X}\Phi,\Phi](Y). \label{CV}
\end{equation}
Its vanishing will tell us that eigenspaces of $\Phi$ are spanned by basic vector fields, i.e.\ vector fields on $M$. Finally, when $\Phi$ cannot tell us much because it is a multiple of the identity, we need assistance of another tensor, called the {\sl tension\/} ${\bf t}$, a type (1,1) tensor field along $\tau$, with components
\begin{equation}
{\bf t}^i_j = \G^i_j - v^k\fpd{\G^i_j}{v^k}. \label{t}
\end{equation}

{\bf Main Theorem}. Assume that $\Phi$ is diagonalizable {\sl with real eigenfunctions} and that
\begin{equation}
\CV{\Phi}=0 \qquad [\nabla\Phi,\Phi]=0\qquad R=0, \label{thm}
\end{equation}
then the system separates into single equations, one for each  1-dimensional eigenspace of $\Phi$ and into individual subsystems for each degenerate eigenfunction, which is then necessarily constant. Each subsystem further decouples if, and only if, the tension $\tt$ is diagonalizable (over $\R$) and $ \CV{\tt}=0$.

\underline{Remarks}:
\begin{itemize}\vspace{-2mm}
\item[-] Constructing appropriate coordinates boils down to integrating distributions which are Frobenius integrable.
\item[-] It is the $\del$-invariance which guarantees that the basic eigendistributions of $\Phi$ are integrable and that connection coefficients relating to different eigendistributions are zero.
\item[-] The zero curvature guarantees that degenerate eigenvalues are constant.
\end{itemize}
Finally, the following general result, also established in \cite{MCS3}, will provide extra information when we turn to the special case of linear systems.

{\bf Proposition}.  Let $D$ be a self-dual derivation. The eigendistributions of a diagonalizable type (1,1) tensor field $U$ along $\tau$ are invariant under $D$ if, and only if, $[DU,U]=0$. It follows that $DU$ is diagonalizable with eigenvalues which are the $D$-derivatives of the eigenvalues of $U$.

So let's have a look at the particular case of linear systems now, when the equations (\ref{nlsode}) get the form (\ref{1}):
\begin{equation}
\ddot{q}^i = - C^i_j \dot{q}^j - K^i_j q^j. \label{lsode}
\end{equation}
It is straightforward to verify from the above coordinate expressions that the tensorial objects figuring in the main theorem take the form
\begin{equation}
   \Phi = K - \onefourth C^2, \qquad \del \Phi = \onehalf [C,K],  \qquad \tt = \left( \G^i_j \right) = \onehalf C, \label{ls1}
\end{equation}
and obviously
\begin{equation}
R\equiv 0, \qquad \CV{\Phi} \equiv 0. \qquad \CV{\tt} \equiv 0. \label{ls2}
\end{equation}
But the above proposition gives us some extra information. Indeed, with $\del$ and $\Phi$ in the role of $D$ and $U$, since the eigenvalues of $\Phi$ will be constants and thus have zero $\del$-derivatives, it tells us that a condition like $[\del\Phi,\Phi] = 0$ will actually imply that $\del \Phi=0$.

The conclusion for linear systems therefore reads: If $K-\onefourth C^2$ has distinct real eigenvalues, then a coordinate transformation will decouple the \sode, provided that $[C,K]=0$; if some eigenvalues are degenerate, but $K-\onefourth C^2$ remains diagonalizable, then decoupling can be achieved, provided that also $C$ is diagonalizable. Obviously, this is consistent with the case of {\sl modal analysis\/} for decoupling.

Just a few words about the extensions in \cite{CSVM} and \cite{ST}. The extension in \cite{CSVM} is in a way the fully time-dependent version of what we have sketched above. But it has an impact on autonomous equations as well, in the sense that there are examples of autonomous systems which fail to satisfy the above requirements for complete separability, but can be decoupled anyway if time-dependent transformations are allowed into the picture. On the other hand, example 2 of the preceding section is a case where $\Phi$ is a multiple of the identity, but the tension $\tt$ happens to have complex eigenvalues so that the above result does not apply. The extension in \cite{ST} covers the case of complex eigenfunctions for $\Phi$ or $\tt$ and leads to a conclusion where `separability' does not refer to single equations, but to pairs of equations which together form a single complex equation. For the equations (\ref{ex2}), it means that after a suitable change of coordinates, the system takes the form of the single complex equation
\begin{equation}
\ddot{z} = -(1+i)z + (1+2\,i)\dot{z}. \label{complex2}
\end{equation}
The technique of {\sl phase synchronization\/} of course has shown us that complete decoupling can be achieved anyway, in some sense. We now finally develop a model that explains how one can interpret this in geometrical terms. For the current situation, where we are dealing with linear differential equations and a linear transformation such as (\ref{tf}), it would be sufficient to talk about \sode s living on $\R^{2n}$ and to see (\ref{tf}) as a coordinate change on that space. But we choose to describe the process in a more general set-up, just in case there would be a chance for generalizations in the future. Roughly, the idea is that we are looking at two different \sode s, potentially living on different tangent bundles, and then a diffeomorphism $f$ between those manifolds which makes that the two \sode\ vector fields are, what is called, $f$-related. This means that the scheme drawn below is commutative, in the sense that $Tf\circ \Gamma_D = \Gamma_C \circ f$.

\setlength{\unitlength}{.8cm}
\begin{picture}(5,5)(-6,0)
\put(3.75,1){\vector(-1,0){2.5}} \put(3.75,4){\vector(-1,0){2.5}}
\put(.5,1.6){\vector(0,1){1.8}}
\put(4.5,1.6){\vector(0,1){1.8}}
\put(2.5,.6){\makebox(0,0){$f$}}
\put(.5,1){\makebox(0,0){$TM$}} \put(.5,4){\makebox(0,0){$TTM$}}
\put(4.5,1){\makebox(0,0){$TN$}} \put(4.5,4){\makebox(0,0){$TTN$}}
\put(5,2.5){\makebox(0,0){$\Gamma_D$}}
\put(0.1,2.5){\makebox(0,0){$\Gamma_C$}}
\put(2.6,4.35){\makebox(0,0){$Tf$}}
\end{picture}

Let's denote coordinates on the tangent bundle $TM$ by $(q^i,\dot{q}^i)$ now; likewise, coordinates on $TN$ are $(p^i,\dot{p}^i)$ say; $f$ represents the map (\ref{tf}), and $Tf$ is its associated tangent map; $\Gamma_C$ is a notation for the \sode\ (\ref{1}) with damping matrix $C$, and $\Gamma_D$ likewise refers to (\ref{eqp}). As is well known, an equivalent way of expressing that the vector fields are $f$-related is saying that for any function $h(q,\dot{q})$ on $TM$, we have that $\Gamma_C(h)\circ f = \Gamma_D(h\circ f)$. Let's see what that tells us in coordinates. We have
\[
\Gamma_C(h) = \dot{q}^i \fpd{h}{q^i} - (C^i_j\dot{q}^j + K^i_j q^j) \fpd{h}{\dot{q}^i},
\]
and subsequently
\begin{equation}
\Gamma_C(h)\circ f = ({T_3}^i_k p^k + {T_4}^i_k \dot{p}^k) \fpd{h}{q^i}
  -  \Big( C^i_j({T_3}^j_k p^k + {T_4}^j_k \dot{p}^k)+ K^i_j ({T_1}^j_k p^k + {T_2}^j_k \dot{p}^k)\Big) \fpd{h}{\dot{q}^i}, \label{lhs}
\end{equation}
where the arguments in the derivatives of $h$ should be thought of now as being expressed in terms of the $(p,\dot{p})$ also. Likewise, we have
\begin{equation}
\Gamma_D(h\circ f) = \dot{p}^k \Big(\fpd{h}{q^i} {T_1}^i_k + \fpd{h}{\dot{q}^i} {T_3}^i_k\Big) - (D^j_k \dot{p}^k + B^j_k p^k) \Big(\fpd{h}{q^i} {T_2}^i_j + \fpd{h}{\dot{q}^i} {T_4}^i_j\Big). \label{rhs}
\end{equation}
The function $h$ being arbitrary, we have to identify the coefficients of $\partial h/\partial q^i$ in (\ref{lhs}) and (\ref{rhs}), and likewise for the coefficients of $\partial h/\partial \dot{q}^i$. This produces two linear expressions in $p^k$ and $\dot{p}^k$, which have to be identically satisfied. In turn, the implication is that the respective coefficients have to be zero. It is straightforward to verify that this gives rise to four matrix relations, which are exactly the conditions (\ref{id1}) and (\ref{id2}) derived in Section~3.

\section{Concluding remarks}

The technique of phase synchronization appears to be an effective computational tool to establish a link between an arbitrary system of linear second-order differential equations and a related system of decoupled equations. It is not a surprise that this technique does not make its appearance in earlier work on separability of general (non-linear) second-order equations, because the differential geometric framework for those studies does not allow transformations depending on velocities. However, from the purely analysis point of view, nobody seems to have picked up certain aspects of ambiguity in the phase synchronization analysis. To begin with, there is the question: how should one understand that the linear relations (\ref{tf}) ``transform'' the given equations (\ref{1}) into the decoupled equations (\ref{eqp})? The point there is that one can think of two systems of second-order equations coming from manipulations with (\ref{tf}), one which has nothing to do with the given equations and one which does, so there is an issue of compatibility. Another question is that there seems to be a lot of freedom in selecting matrices $T_i$ which have the same effect. We have illustrated various aspects of these questions. It seems to us that there is room for a deeper study of the compatibility question in the format (\ref{id1},\ref{id2}), which could lead for example to an understanding of the mathematical structure of the set of $T_i$ having that same effect. And since we managed to give a geometrical interpretation of these conditions, a more challenging idea is to look for possible generalizations to non-linear equations. Allow us finally a more critical remark: we have gratefully used survey formulas from \cite{SKML18} for notational convenience, but it seems to us, for example in view of the non-uniqueness of the $T_i$, that it is highly questionable whether linking a transformation like (\ref{tf}) to a kind of generalized inverse problem of Lagrangian mechanics can make any sense!

\section*{Acknowledgments}

TM thanks the Research Foundation -- Flanders (FWO) for its support through Research Grant 1510818N.




\begin{thebibliography}{99}

\bibitem{Bloch}
A.M.\ Bloch, P.S.\ Krishnaprasad, J.E.\ Marsden and T.S.\ Ratiu,
Dissipation induced instabilities,
{\em Ann.\ Inst.\ H.\ Poincar\'e. Analyse Nonlineare} {\bf 11} (1994) 37--90.


\bibitem{CSVM}
F.\ Cantrijn, W.\ Sarlet, A.\ Vandecasteele and E.\ Mart\'{\i}nez,
Complete separability of time-dependent second-order ordinary differential equations,
{\em Acta Applic.\ Math.\/} {\bf 42} (1996) 309--334.

\bibitem{Goldstein}
H.\ Goldstein, {\it Classical Mechanics\/} (2nd.\ edition)
(Addison-Wesley 1980).


\bibitem{KMM11}
T.\ Kawano, M.\ Morzfeld and F.\ Ma, The decoupling of defective linear dynamical systems in free motion, {\em J.\ Sound Vib.\/} {\bf 330} (2011) 5165--5183.

\bibitem{MIM09}
F.\ Ma, A.\ Imam and M.\ Morzfeld, The decoupling of damped linear systems in oscillatory free
vibration, {\em J.\ Sound Vib.\/} {\bf 324} (2009) 408--428.

\bibitem{MMI10}
F.\ Ma, M.\ Morzfeld and A.\ Imam, The decoupling of damped linear systems in free or forced
vibration, {\em J.\ Sound Vib.\/} {\bf 329} (2010) 3182--3202.

\bibitem{MCS1}
E.\ Mart\'{\i}nez, J.\ F.\ Cari\~{n}ena and W.\ Sarlet,
Derivations of differential forms along the tangent bundle projection,
{\em Diff.\ Geometry and its Applications\/} {\bf 2} (1992) 17--43.

\bibitem{MCS2}
E.\ Mart\'{\i}nez, J.\ F.\ Cari\~{n}ena and W.\ Sarlet,
Derivations of differential forms along the tangent bundle projection. Part II,
{\em Diff.\ Geometry and its Applications\/} {\bf 3} (1993) 1--29.

\bibitem{MCS3}
E.\ Mart\'{i}nez, J.\ F.\ Cari\~{n}ena and W.\ Sarlet,
Geometric characterization of separable second-order differential equations,
{\em Math.\ Proc.\ Camb.\ Phil.\ Soc.\/} {\bf 113} (1993) 205--224.

\bibitem{MM11}
M.\ Morzfeld, F.\ Ma, The decoupling of damped linear systems in configuration
and state spaces, {\em J.\ Sound Vib.\/} {\bf 330} (2011) 155--161.

\bibitem{MMP11}
M.\ Morzfeld, F.\ Ma and B.N.\ Parlett, The transformation of second-order linear systems into
independent equations, {\em SIAM J.\ Appl.\ Math.\/} {\bf 71} (2011) 1026--1043.


\bibitem{Rosenberg}
R.\ M.\ Rosenberg, {\it Analytical Dynamics of Discrete Systems\/}
(Plenum Press 1977).


\bibitem{SKML18}
R.G.\ Salsa Jr., D.T.\ Kawano, F.\ Ma and G.\ Leitmann, The inverse problem of linear Lagrangian dynamics, {\em ASME J.\ Appl.\ Mech.\/} {\bf 85} (2018) 031002.


\bibitem{Sreview} W.\ Sarlet, Complete decoupling of systems of ordinary second-order differential equations, In: N.H.\ Ibragimov, F.M.\ Mahomed, D.P.\ Mason and D.\ Sherwell (Eds.), Proc.\ 4th Workshop on Differential Equations and Chaos (1997) 237--264.

\bibitem{S}
W.\ Sarlet, Different forms of separability of second-order equations, {\em Nonlinear Analysis\/} {\bf 47} (2001) 6135--6146.


\bibitem{ST}
W.\ Sarlet and G.\ Thompson,
Complex second-order differential equations and separability,
{\em Applicable Algebra in Engineering, Communication and Computing\/} {\bf 11} (2001) 333--357.



\end{thebibliography}
\end{document}